\documentclass{lmcs}


\usepackage{hyperref}

\newcommand{\Prob}{{\mbox{Prob}}}
\newcommand{\bits}{{\{0,1\}}}
\newcommand{\uu}{{\bits^*}}
\newcommand{\nn}{{\{0,1\}^n}}
\newcommand{\tr}{{|\mkern-2.5mu|}}
\newcommand{\np}{\mbox{{NP}}}
\newcommand{\npco}{{{\np \cap \conp}}} % chktex 1
\newcommand{\conp}{{\mbox{coNP}}}

\usepackage[utf8]{inputenc}

\newtheorem*{rems*}{Remarks}
\newtheorem{myclm}{Claim}

\usepackage{lastpage}
\lmcsdoi{16}{3}{9}
\lmcsheading{}{\pageref{LastPage}}{}{}%
{Apr.~07,~2020}{Aug.~12,~2020}{}

\begin{document}

\title{A limitation on the KPT interpolation}

\author{Jan Kraj\'{\i}\v{c}ek}

\address{MFF, Charles University, Sokolovsk\' a 83, Prague, 186 75, The Czech Republic}

\email{krajicek@karlin.mff.cuni.cz}

\begin{abstract}
\noindent
We prove a limitation on a variant of the KPT theorem proposed for propositional proof systems
by Pich and Santhanam~\cite{PichSan20}, for all proof systems that prove the disjointness of two
NP sets that are hard to distinguish.
\end{abstract}

\maketitle

For a $\conp$ property $\psi(x)$, given $n \geq 1$, we can construct a size $n^{O(1)}$ propositional formula $\tr \psi \tr^n(x,y)$
with $n$ atoms $x = (x_1, \dots, x_n)$ and $n^{O(1)}$ atoms $y$ such that for any $a \in \nn$, $\psi(a)$ is true iff
$\tr \psi\tr^n(a,y) \in \mbox{TAUT}$. This is just a restatement of the $\np$-completeness of SAT\@. In addition, if $\psi(x)$ is defined
in a suitable language of arithmetic and has a suitable logical form, the translation can be defined purely syntactically without
a reference to machines or computations. This then allows to transform also a possible first-order proof of $\forall x \psi(x)$
into a sequence of short propositional proofs of tautologies $\tr \psi\tr^n$, $n = 1, 2, \dots$; if the original proof uses axioms of theory $T$
(essentially any sound r.e.\ theory) then the propositional proofs will be in a proof system $P_T$ associated to $T$. Many standard proof systems are of the form $P_T$ for some $T$, and this is often
the most efficient way how to construct short $P_T$-proofs of uniform sequences of tautologies.
Although the unprovability of $\forall x \psi(x)$ in $T$ does not imply lower bounds for $P_T$-proofs
of the tautologies, a method used in establishing the unprovability sometimes yields an insight how the lower bound could be proved.
All this is a well-established part of proof complexity and the reader can find it in~\cite[Chpt.12]{prf} (or in references given there).

The translation is, however, not entirely faithful for formulas of a certain logical form, and this is an obstacle for transforming
the conditional unprovability result for strong universal
theories in~\cite{Kra-nwg} into conditional lower bounds for strong proof systems.
To explain the problem in some detail assume $\psi(x)$ has the
form
\begin{equation} \label{15.3.20a}
\exists i < |x| \forall y (|y| = |x|)\ \varphi(x,i,y)
\end{equation}
where $\varphi$ is a p-time property and $|x|$ is the bit length of $x$.
The provability of $\forall x \psi(x)$ in a universal $T$ can be analyzed using the
KPT theorem which provides an efficient interactive algorithm for finding $i$ given $x$
(cf.~\cite{KPT} or~\cite[Sec.12.2]{prf}). The same method does not, however, work in the propositional setting. To illustrate this assume that $\tr \psi\tr^n$ has a proof in proof system $P_T$ attached to $T$ and from that we can deduce in $T$ that
\begin{equation} \label{15.3.20b}
\bigvee_{i < n}\ \tr\psi\tr^n(x,i,y_i)
\end{equation}
is a tautology (in addition the translation assures that all $y_i$ are disjoint tuples of atoms). This
implies in $T$ that for all assignments $a$ and $b = (b_0, \dots, b_{n-1})$ for all $x$ and all $y$ variables there is $i < n$ such that
$\tr \psi\tr^n(a,b)$ is true. But to get (\ref{15.3.20a}) (and then use the KPT analysis from~\cite{Kra-nwg})
we would need to show that
for all $a$ there is one $i < n$ such that for all $b_i$ the formula is true. Unfortunately, to derive this one needs to use
the bounded collection scheme (allowing to move the quantifier bounding $i$ before the quantifier bounding $b_i$)
and this scheme is not available
in universal theories under consideration, cf.~\cite{CT}. The reader can find more about this issue in~\cite[Sec.5]{Kra-nwg}
or at the end of~\cite[Sec.12.8]{prf};
knowing this background offers my motivation for this research (which differs perhaps from that of~\cite{PichSan20})
but it is not needed to understand the argument below.

Pich and Santhanam~\cite{PichSan20} proposed a direct way how to bypass this obstacle: simply ignore it and prove a version of
the KPT theorem for (some, at least) strong propositional proof systems. For such proof systems a conditional
lower bound can be indeed proved, cf.~\cite{PichSan20} or~\cite{Kra-nwg}.

\begin{defiC} [\cite{PichSan20}]
Let P be a propositional proof system. The system has {\bf KPT interpolation} if there are a constant $k \geq 1$
and $k$ p-time functions
\[
f_1(x,z), f_2(x,z,w_1), \dots, f_k(x, z, w_1, \dots, w_{k-1})
\]
such that whenever $\pi$ is a P-proof of a disjunction of the form
\[
A_0(x,y_1) \vee \dots \vee A_{m-1}(x,y_m)
\]
where $x$ is a $n$-tuple of atoms and $y_1, \dots, y_m$ are disjoint tuples of atoms, then for all
$a \in \nn$ the following is valid for all $b_1, \dots, b_m$ of the appropriate lengths:
\begin{itemize}

\item either $A_{i_1}(a,y_{i_1}) \in \mbox{TAUT}$ for $i_1 = f_1(a,\pi)$ or, if $A_{i_1}(a,b_{i_1})$ is false,

\item $A_{i_2}(a,y_{i_2}) \in \mbox{TAUT}$ for $i_2 = f_{2}(a, \pi, b_{i_1})$ or, if $A_{i_2}(a,b_{i_2})$ is false,

\item $\dots$, or

\item $A_{i_k}(a,y_{i_k}) \in \mbox{TAUT}$ for $i_k = f_{k}(a, \pi, b_{i_1}, \dots, b_{i_{k-1}})$.
\end{itemize}
\end{defiC}

\noindent
An illuminating interpretation of the definition can be made using the interactive communication model of~\cite{KPS} involving Student and Teacher. Student is a p-time machine while Teacher has unlimited powers.
At the beginning Student gets $a \in \nn$ and the proof $\pi$ and computes from it his first candidate solution:
index $i_1$ such that $A_{i_1}(a, y_{i_1})$ is --- he thinks --- a tautology. Teacher either approves or she
provides Student with a counter-example: an assignment $b_{i_1}$ for $y_{i_1}$ which falsifies the formula.
In the next round Student can use this counter-example to propose his next candidate solution, etc. Functions
$f_1, \dots f_k$ in the definition form a strategy for Student so that he solves the task for all $a$ and $\pi$ in
$k$ steps in the worst case. Note that if we fixed $m=2$ as in ordinary interpolation then $k=2$ would suffice;
the concept makes sense for variable $m$ only.

\bigskip

Unfortunately, we show in this note that this property fails for strong proof systems (above a low depth Frege system)
for essentially the same reasons why ordinary feasible interpolation fails for them (cf.~\cite[Sec.18.7]{prf}).
For a set $U \subseteq \uu$ and $n \geq 1$ put $U_n := U \cap \bits^n$. $\mbox{LK}_{3/2}$ is the $\Sigma$-depth $1$ subsystem of
sequent calculus (cf.~\cite[Sec.3.4]{prf}).

\begin{thm}
Let P be a proof system containing $\mbox{LK}_{3/2}$.
Assume that $U, V$ are disjoint NP sets such that:
\begin{enumerate}

\item Propositional formulas expressing that $U_n \cap V_n = \emptyset$ have p-size P-proofs.

\item For any constant $c \geq 1$, for all large enough $n$ there is a distribution ${\bf D}_n$
on $\nn$ with support $U_n \cup V_n$ such that there is no size $n^c$ circuit $C_n$
for which
\[
\Prob_x[(x \in U_n \wedge C_n(x) = 1) \vee (x \in V_n \wedge C_n(x) = 0)] \geq 1/2 + n^{-c}
\]
where samples $x$ in the probability are chosen according to ${\bf D}_n$.

\end{enumerate}

\noindent
Then P does not admit KPT interpolation.

\end{thm}

\begin{rems*}
\hfill
\begin{enumerate}
\item An example of a pair of two NP sets $U,V$ that are conjectured to satisfy the second
condition can be defined using one-way permutation (more generally an injective one-way function with
output length determined by input length)
and its hard bit: $U$ (resp. $V$) are the strings in the range of the permutation
whose hard bit is $1$ (resp. $0$). Distribution ${\bf D}_n$ is in this case generated by the permutation
from the uniform distribution on the seed strings, i.e.\ it is uniform itself.

\item It is known that the hypothesis of the theorem can be fulfilled for systems such as EF, F, $\mbox{TC}^0$-F and,
under stronger hypotheses about non-separability of $U$ and $V$, also for $\mbox{AC}^0$-F above certain small depth; see
the comprehensive discussion in~\cite[Sec.18.7]{prf}.

\item The phrase that P contains $\mbox{LK}_{3/2}$ means for simplicity just that: P can operate with sequents
consisting of $\Sigma$-depth $1$ formulas and all $\mbox{LK}_{3/2}$-proofs are also P-proofs.
However, this is used only in Claim~\ref{claim:first} and, in fact,
it would suffice that P represents formulas $U(x,y)$ and $V(x,z)$ (defined below) in some other formalism and efficiently simulates modus ponens.

\end{enumerate}
\end{rems*}

\bigskip
\noindent
{\bf Proof of the theorem} occupies the rest of this note.

\medskip

Write $U(x,y)$ for a p-time relation that $y$ witnesses $x \in U$
and similarly $V(x,z)$ for $V$, with the length of both $y$ and $z$ p-bounded in the length of $x$.
Let $n,m \geq 1$ and for $m$ strings $x_1, \dots x_m$ of length $n$ each
consider the following $2m$ propositional formulas translating the predicates $U(x,y)$ and $V(x,z)$
(which we shall denote also $U$ and $V$ in order to ease on notation):
\begin{itemize}

\item $U(x_i, y_i)$: $x_i$ is an $n$-tuple of atoms for bits of $x_i$ and
$y_i$ is an $n^{O(1)}$-tuple of atoms for bits of a witness associated with $x_i$ together with bits needed to
encode $U$ as propositional formula suitable for P (e.g.\ as 3CNF),

\item $V(x_i, z_i)$: analogously for $V$,

\item where all $x_i, y_i, z_i$ are disjoint.
\end{itemize}
Consider the induction statement:
\begin{equation} \label{3a}
 x_1 \in U \wedge (\forall i < m,\ x_i \in U \rightarrow x_{i+1} \in U)
 \rightarrow x_m \in U
\end{equation}
and write it as a disjunction with $m+1$ disjuncts:
\begin{equation}\label{3b}
 x_1 \notin U \vee
 \bigvee_i (x_i \in U \wedge x_{i+1} \notin U)
 \vee x_m \in U\ .
 \end{equation}
Now replace $x_i \in U$ by $x_i \notin V$
and $x_m \in U$ by $x_m \notin V$ and write it propositionally:
\begin{equation} \label{1}
   \neg U(x_1, y_1) \vee
   \bigvee_i [\neg V(x_i , z_i) \wedge \neg U(x_{i+1}, y_{i+1})]
   \vee  \neg V(x_m, z_m).
\end{equation}
Note that except the $x$-variables the $m+1$ disjuncts are disjoint.

\begin{myclm}%
\label{claim:first}
(\ref{1}) has a p-size proof in P.
\end{myclm}

\smallskip
\noindent
To see this note that
induction (\ref{3a}) can be proved by simulating modus ponens (here we use that P contains $\mbox{LK}_{3/2}$).
Disjunction (\ref{1}) follows from it because we assume that the disjointness of $U_n, V_n$ has short P-proofs,
i.e. $U(x,y) \rightarrow \neg V(x,z)$ has a short proof.

\bigskip

Now apply the supposed KPT interpolation to (\ref{1}). W.l.o.g.
we shall assume (and arrange that in the construction below)
that $x_1 \in U$ and $x_m \in V$ (with witnesses $y_1$ and $z_m$, respectively).
Hence Student in the KPT computation is
supposed to find $i < m$ for which the $i$-th disjunct
\[
A_i\ :=\ [\neg V(x_i , z_i) \wedge \neg U(x_{i+1}, y_{i+1})]\ ,\ i = 1, \dots, m-1
\]
is valid (i.e.\ where the induction step going from $i$ to $i+1$ fails).
We shall show that the existence of such a KPT p-time Student
allows to separate $U_n$ from $V_n$ with a non-negligible advantage violating the hypotheses of the theorem.

Take any $m$ such that $3\cdot 2^{k-1} \le m \le n^{O(1)}$ (the upper bound implies that the proof in Claim~\ref{claim:first} is of size
$n^{O(1)}$). For $1 \le i < m$ define:
\[
W_i[m] := U^i \times V^{m-i}\ \mbox{ and } \
W[m] := \bigcup_i W_i[m]\ .
\]
Note that any string $w = (w_1, \dots, w_m) \in W[m]$ satisfies $w_1 \in U$ and $w_m \in V$.

Let $k \geq 1$ and $f_1, \dots, f_k$ be the constant and the p-time functions provided the assumed KPT interpolation for P.
Assume that $1 \le i_1 < m$ is the most frequent value $f_1$ computes on inputs from $W[m]$ (thinking of a P-proof $\pi$ as
fixed). This maximal frequency $\gamma$
is at least $1/m$. (Here the \emph{frequency} means with respect to the product of distributions ${\bf D}_n$ on $\bits^n$ for which
it is assumed that $U_n, V_n$ are hard to separate.)

\begin{myclm}%
\label{claim:second}
The frequency on $W_{i_1}[m]$ is at least $\gamma - n^{\omega(1)}$, i.e.\ it is at least $1/m$
modulo a negligible error.
\end{myclm}

\smallskip
\noindent
Note that for any $i < j$ the frequency for $W_i[m], W_j[m]$ can differ only negligibly because otherwise we could
use the usual triangle inequality argument to find a non-negligible discrepancy between frequencies on
$W_t[m]$ and $W_{t+1}[m]$ for some $i \le t < j$, and use it to separate $U_n$ from $V_n$
(on position $t+1$, after fixing the rest of coordinates by averaging).
Because all $W_i[m]$ are disjoint, the frequency must be $\gamma$ up to a negligible difference.

\bigskip

Now we describe a process that transforms the assumed successful strategy\! for\! Student into a p-time algorithm\! with\! p-size advice, separating $U_n, V_n$ with a non-negligible advantage.

Assume first $i_1 < m/2$. By averaging there are $u_1, \dots u_{m/2} \in U_n$ s.t. $f_1(w) = i_1$
with frequency at least $1/(2m)$ (the factor $2$ in the denominator allows us to forget about the \emph{``up to the
negligible error''} phrase) for all $w$ of the form:
\[
\{u_1\} \times \cdots \times \{u_{m/2}\} \times W[m/2]\ .
\]
Fix such $u_1, \dots, u_{m/2}$ and also witnesses $a_1, \dots, a_{m/2}$ for their membership in $U$.
These will be used as advice for the eventual algorithm.

If $i_1 \geq m/2$ then fill analogously the last $m/2$ positions by elements of $V_n$ and include the relevant
witnesses in the advice. W.l.o.g.\ we assume that the first case $i_1 < m/2$ occurred.

\smallskip

We interpret this situation as reducing the Student-Teacher computation to $k-1$ rounds on smaller
universe $W[m/2]$. Namely, given $w = (w_1, \dots, w_{m/2}) \in W[m/2]$ define:
\begin{equation} \label{22.3.20a}
\tilde w\ :=\
(u_1, \dots, u_{m/2}, w_1, \dots, w_{m/2})\ \in \ W[m]
\end{equation}
and run $f_1$ on $\tilde w$. If $f_1(\tilde w) \neq i_1$, declare failure. Otherwise use
the advice witnesses to produce a falsifying assignment for $A_{i_1}$: $U(u_{i_1 + 1}, a_{i_1 + 1})$ holds.

After this first step use functions $f_2, f_3, \dots$
(and Claim~\ref{claim:second} for the smaller universes) and as long as they give values $j < m/2$
always answer for Teacher using the advice strings $a_j$.
Eventually Student proposes value $j \geq m/2$: choose the most frequent such value
$i_2 \geq m/2$ and proceed as in case of $i_1$, further restricting domain (\ref{22.3.20a}) as in
binary search.
Repeating this at most $(k-1)$-times the situation will be as follows:
\begin{enumerate}

\item The universe will shrink at most
to $W[m/(2^{k-1})]$ which is at least $W[3]$. In fact, we shall arrange in the last step
that exactly $W[3]$ remains (by filling in more positions by elements of $U_n$ or $V_n$, respectively, if needed)
and hence the inputs before applying the last KPT function $f_k$
are of the form $(w_1, w_2, w_3)$ with $w_1 \in U$ and $w_3 \in V$.

Note that Student gets to use $f_k$ because if he succeeded earlier it would violate Claim~\ref{claim:second}.

\item The last function $f_k$
has to find a gap in the induction, and this itself will violate Claim~\ref{claim:second}.
In particular, the gap is either between $w_1$ and $w_2$ and then $w_2 \in V$,
or between $w_2$ and $w_3$ and then $w_2 \in U$.

\item This process has the probability $\geq 1/(2m)$, i.e.\ non-negligible, of not failing in any of the $k-1$ rounds and hence
it will not fail and will compute correctly the membership of (any) $w_2$ in $U$ or $V$ with a non-negligible probability.
In all cases when the process fails output random bit $0$ or $1$ with equal probability.

\end{enumerate}
This proves the theorem.
\qed%

\bigskip

\noindent
We conclude by pointing out that the KPT theorem enters propositional proof complexity also via notions of
pseudo-surjective and iterable maps in the theory of proof complexity generators, cf.~\cite{Kra-dual}
or~\cite[Sec.19.4]{prf} for detailed expositions of this subject.

\section*{Acknowledgment}
  \noindent I thank J.~Pich (Oxford) for comments on an earlier note.

\end{document}